\documentclass[12pt]{article}
\usepackage{graphics}
\usepackage{color}
\usepackage{epsfig}
\usepackage{latexsym}
\usepackage{amssymb}
\usepackage{multicol}
\usepackage{graphicx}
\usepackage{amsmath}
\usepackage{amsfonts}
\usepackage{enumerate}
\DeclareMathSymbol\square {\mathrel}{AMSa}{"03}
\begin{document}
\title {Pythagorean Triples and A New Pythagorean Theorem}
\author{H. Lee Price and Frank Bernhart}
\date{January 1, 2007}


\maketitle

\begin{abstract}
\textit{Given a right triangle and two inscribed squares, we show that the 
reciprocals of the hypotenuse and the sides of the squares satisfy an 
interesting Pythagorean equality. This gives new ways to obtain rational 
(integer) right triangles from a given one.}  
\end{abstract}

\begin{flushleft}
\textbf{1. Harmonic and Symphonic Squares}
\end{flushleft}

Consider an arbitrary triangle with altitude $\alpha $ corresponding to base 
$\beta $ (see \emph{Figure 1a}). Assuming that the base angles are 
acute, suppose that a square of side $\eta $ is inscribed as shown in \emph{Figure 1b}.

\begin{figure}[h]
\centerline{\includegraphics[width=2.88in,]{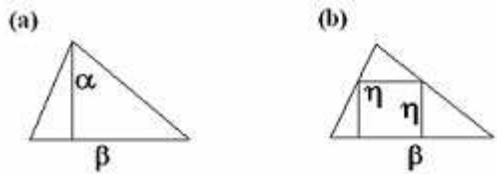}}
\label{fig1}
\caption{\emph{Triangle and inscribed square}}
\end{figure}

Then $\alpha ,\beta ,\eta $ form a \emph{harmonic sum}, i.e. satisfy \textit{(\ref{eq1})}. The equivalent formula $\eta =\frac{\alpha \beta }{\alpha +\beta }$ 
is also convenient, and is found in some geometry books.
\begin{equation}
\label{eq1}
\frac{1}{\alpha }+\frac{1}{\beta }=\frac{1}{\eta }
\end{equation}
Equation \textit{(\ref{eq1})} remains valid if a base angle is a right angle, or is obtuse, save that in 
the last case the triangle base must be extended, and the square is not 
strictly inscribed ( \emph{Figures 2a, 2c} ).

\begin{figure}[h]
\centerline{\includegraphics[width=3.825in,]{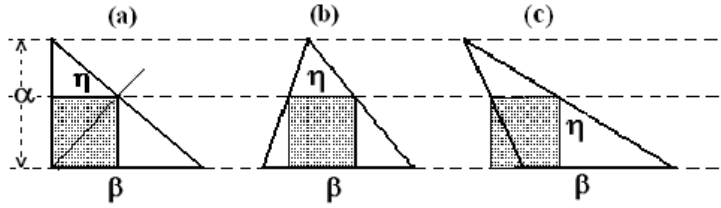}}
\label{fig2}
\caption{ Three congruent harmonic squares}
\end{figure}

Starting with the right angle case ({\it Figure 2a}) it is easy to see the inscribed square 
uniquely exists (bisect the right angle), and from similar triangles we have 
proportion
\[
\left( {\beta -\eta } \right):\beta =\eta :\alpha 
\]
which leads to \textit{(\ref{eq1})}. The horizontal dashed lines are parallel, hence the three 
triangles of {\it Figure 2} have the same base and altitude; and the three squares are 
congruent and unique. Clearly \textit{(\ref{eq1})} applies to all cases.

A scalene triangle will have three squares, one for each side of the 
triangle. Each one has a side length that is the harmonic sum of the 
corresponding triangle side and altitude. Of course when the triangle is 
equilateral, the squares are congruent, but do not coincide.

When a right triangle with legs $a,b$ and hypotenuse $c$  is given, 
there are just two squares ({\it Figures 3a, 3b}), the ``harmonic'' square of side $h$, and the 
``symphonic'' square of side $s$. 

\begin{figure}[h]
\centerline{\includegraphics[width=3.00in]{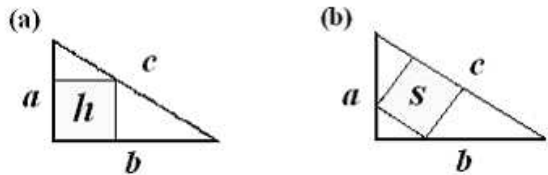}}
\label{fig3}
\caption{\emph{Harmonic} and \emph{symphonic} squares.}
\end{figure}

Using \textit{(\ref{eq1})} twice gives solutions \textit{(\ref{eq2})}, since the altitude to the hypotenuse is 
$\kappa=\frac{ab}{c}$(using similar right triangles).
\begin{equation}
\label{eq2}
h=\frac{ab}{a+b},\quad s=\frac{abc}{ab+c^2}
\end{equation}
In all right triangles (Exercise) $c>h>s$.

Now we present a striking identity we call  \emph{the Symphonic Theorem}.

\begin{flushleft}
{\bf The Symphonic Theorem:}
\end{flushleft}
\textit{With respect to $c$, $h$, and $s$ as defined above, the triple $[c^2,h^2,s^2]$ is harmonic, and the triple
$[ \tfrac{1}{c},\tfrac{1}{h},\tfrac{1}{s}]$ is Pythagorean.}
\begin{equation}
\label{eq3}
 \frac{1}{c^2}+\frac{1}{h^2}=\frac{1}{s^2}
\end{equation}

Equation \textit{(\ref{eq1})} holds for all right triangles $[a,b,c]$ and squares $h,s$ inscribed as indicated. The altitudes of the triangle (${a}, {b}, {\kappa}$), have a similar relationship. 

\begin{equation}
\label{eq4} 
\frac{1}{{a}^2}+\frac{1}{{b}^2}=\frac{1}{{\kappa }^2}
\end{equation}
Substitute \textit{(\ref{eq2})} into \textit{(\ref{eq3})} and clear fractions:
\begin{equation}
\label{eq5}
[c(a+b)]^2+[ab]^2=[ab+c^2]^2.
\end{equation}
On both sides of the equal sign, make the replacement
\begin{equation}
\label{eq6}
c^2=a^2+b^2
\end{equation}
to get an algebraic identity:
\[
[a^2+ b^2][a+b]^2+[ab]^2 = [ab + a^2 + b^2]^2.
\]
Now suppose that $[a,b,c]$ are integers without common 
factors (a {\it Primitive Pythagorean Triple}, or PPT). From equation \textit{(\ref{eq6})}, $a, b, c$ are 
relatively prime in pairs. Then equation \textit{(\ref{eq5})} shows that 
$[c(a+b),ab,ab+c^2]$ is also a PPT. 
Any common factor divides, $ab+c^2,ab,$ hence also $(ab+c^2)-ab,$ which implies a non-trivial common factor of $a$ or 
$b$ and $c$, contrary to assumption.

Take for example $[a,b,c]=[3,4,5]$ and calculate $h=\tfrac{12}{7},s=\tfrac{60}{37}.$ From \textit{(\ref{eq3})} 
we have a rational right triangle 
$[\tfrac{1}{h},\tfrac{1}{c},\tfrac{1}{s}]=[ \tfrac{7}{12},\frac{1}{5},\tfrac{37}{60}].$
 
 Multiply this by $60$ to arrive at the triple $[35,12,37],$ in accord with \textit{(\ref{eq5})}. On the other hand, we can enlarge $a,\mbox{ }b, \mbox{ }c \mbox{ }$ so that $h, \mbox{ }s$ also 
are integers. Multiply by $(11)(37)= 259,$ to get 
$[a,b,c]=[777,1036,1295],$ $h=444,$ and $s=420.$ This is clearly the smallest example.

If \textit{(\ref{eq3})} is regarded as simply an equation in three unknowns, then 
the smallest solution in integers is $[15,20,12],$ 
derived by dividing the $[3,4,5]$ right triangle by 
a factor of sixty. More generally, given PT $[a,b,c],$ we can obtain $[bc,\mbox{ }ac,\mbox{ }ab]\mbox{ }$as an integer solution to $\frac{1}{x^2}+ \frac{1}{y^2}= \frac{1}{z^2}.$ 

The transformation $\mbox{S: }[3,4,5] \to [35,12,37]$ by \textit{(\ref{eq5})} is an instance of
\begin{equation}
\label{eq7}
\mbox{S: }[a,\mbox{ }b,\mbox{ }c]\mbox{ }\to \mbox{ }[c(a+b),\mbox{ 
}ab,\mbox{ }c^2 + ab].
\end{equation}
Looking at this form one day, we found the following analogous 
transformation.
\begin{equation}
\label{eq8}
\rm S': [a,\mbox{ }b,\mbox{ }c]\mbox{ }\to \mbox{ }[c\vert a-b\vert 
,\mbox{ }ab,\mbox{ }c^2 - ab]
\end{equation}
If in the definitions \textit{(\ref{eq2})} the smaller of the two values 
$a,\mbox{ }b\mbox{ }$ is replaced by its negative, defining $h',\mbox{ }s'$, 
then the same development which led to \textit{(\ref{eq7})} now leads to \textit{(\ref{eq8})}.

The results \textit{(\ref{eq7})} \textit{(\ref{eq8})} are here named the {\bf {\it Symphonic 
Derivatives}}, {\bf {\it Major}} and {\bf {\it Minor}}. One of us (Price) 
found the Theorem, both of us worked out the results, and one of us (Bernhart) 
found the Minor derivative.

A diagram for triangle $[a,b,c]$ may be embellished with lines showing how to find $h'$ and $s'$, but these 
constructions appear highly artificial, compared with the elegant beauty of \textit{Figure 3}! Query: can a ``natural'' construction for \textit{(\ref{eq8})} be devised?

There is one clue, a challenge problem posed by Sastry in \textit{The 
College Mathematics Journal} [10]. The editors made a 
brief composite of five independent solutions, and appended a long list of 
other solvers. The question starts with the two triples $[3,4,5]$ and $[5,12,13],$ and seeks to generalize to 
$[a,b,c]$ and $[c,ab,ab+1].$ This is the special case of the minor derivative when $a-b=1$! This last condition defines a family of primitive triples studied by Fermat (see Eckert [6]).

That family is closely related to a venerable problem of recreational 
mathematics, the square-triangle problem. The challenge is to find pairs of 
positive integers $x,y$ such that a square arrangement of marbles, $x$ on a 
side, equals a triangular arrangement, $y$ on a side.

The question amounts to finding integer solutions to $2x^2=y(y+1),$ which can be converted to a Pellian and solved (see Barbeau [2]). The pairs $(x_i, y_i )$ form 
an infinite family, but we shall focus on 
\[
(x_i )=(1,\mbox{ }6,\mbox{ }35,\mbox{ }204,\mbox{ 
}1189,\cdots,\mbox{ }a,\mbox{ }b,\mbox{ }6b-a,\cdots),
\]
easily extended by the embedded rule $ c=6a-b$. Pick two consecutive values, such as 
$6,35$ or $35,204.$ Take the sum and difference, but 
split the sum into two consecutive integers: 
\[
\begin{array}{rrrr}
 35 - 6 = 29, & \quad 35+6   =  20 + 21   & \mbox{ }\to \mbox{ } & [20,\mbox{ }21,\mbox{ }29]\mbox{ } \\ 
 204- 35=169, & \quad 204+35 = 119 + 120  & \mbox{ }\to \mbox{ } & [119,\mbox{ }120,\mbox{ }169]. \\ 
 \end{array}
\]
We summarize the effects of ``Symphonic derivation'' on a few smaller 
Pythag-orean triangles. 
\[
\begin{array}{lcccc}
 \rm S,S': &   [3,4,5]  & \to  & [35,12,37],    \mbox{ } & [5,12,13]     \\ 
 \rm S,S': & [5,12,13]  & \to  & [221,60,229],  \mbox{ } & [91,60,109]   \\ 
 \rm S,S': & [15,8,17]  & \to  & [391,120,409], \mbox{ } & [119,120,169] \\ 
 \rm S,S': & [7,24,25]  & \to  & [775,168,793], \mbox{ } & [425,168,457] \\ 
 \end{array}
\]
\[
<>\quad <>\quad <>\quad <>\quad <>\quad <>\quad <>\quad <>\quad <>
\]
{\bf 2. Raising the Standard}

It is very helpful to study symphonic derivation from the point of view of 
the standard generators of triples. From one point of view, a pair 
$(p,\mbox{ }q)\mbox{ }$of parameters with certain restrictions is used to 
generate a triple $[a,b,c].$ From another viewpoint, 
a proper fraction $\frac{q}{p}$ in lowest terms, \emph{is} the generator -- with a 
definite geometric meaning. We shall treat the parameter pair and the 
fraction as equivalent, interchangeable objects, with a few 
adjustments. A brief review will be given (details and proofs omitted). A 
fuller treatment is found in [4], [5].

Primitive triple $[a,b,c]$ has two generators. We 
assume without loss that $a$ is odd. {\bf {\it Primary}} and {\bf {\it 
secondary}} generators $t_1 ,\mbox{ }t_{2\mbox{ }} $are obtained as follows.
\[
t_1 \mbox{}= \mbox{}\frac{b}{c+a}\mbox{}=\mbox{}\frac{c-a}{b} \mbox{}=\mbox{}\frac{q_1 }{p_1 },\mbox{  }\quad  \quad \mbox{  }t_2 \mbox{}=\mbox{} \frac{a}{c+b}\mbox{}=\mbox{} \frac{c-b}{a}\mbox{}= \mbox{}\frac{q_2}{p_2}.
\]
Quotients $t_1 ,\mbox{ }t_{2\mbox{ }} $ are geometrically the {\bf {\it 
half-angle tangents}} of the right triangle. The generators for $[3,4,5]$ are: $$t_1=\frac{4}{3+5}=\frac{1}{2},\mbox{  } \quad \quad \mbox{  } t_2=\frac{3}{4+5}=\frac{1}{3}.$$ 

The {\bf {\it key sequence}} $[q_2,q_1,p_1,p_2]$ is a Fibonacci-Rule sequence, $\ q_2+q_1=p_1,$ $q_1+p_1=p_2$ with the additional conditions that the first member $(q_2)$ is odd, and the first two members $(q_1,  q_2)$ are relatively 
prime (and of course all four members are positive integers). Some 
possibilities are as follows.
\[
[1,n,\cdots]\mbox{ }[3,1,\cdots]\mbox{ }[3,2,\cdots]\mbox{ }[3,4,\cdots]\mbox{ }[5,1,\cdots]\mbox{ }[5,2,\cdots] = [5,\mbox{ }2,\mbox{ }7,\mbox{ }9]
\]
We can also make a key sequence from any positive proper fraction 
$\frac{q}{p}.$ If $q+p\mbox{ }$is even, we use the ``template'' $[q,\mbox{ 
}\ast ,\mbox{ }\ast ,\mbox{ }p],\mbox{ }$ but if $q+p\mbox{ }$ is odd, 
we use the template $[\ast ,\mbox{ }q,\mbox{ }p,\mbox{ }\ast ].\mbox{ }${\bf 
} In either case, the template can be completed uniquely. For example, fractions $\frac{2}{3},$ $\frac{1}{5}$ each give the same sequence, $[1,\mbox{ }2,\mbox{ }3,\mbox{ }5],\mbox{ }$ and fractions $\frac{1}{4},$ $\frac{3}{5}^{ }$ each 
give $[3,\mbox{ }1,\mbox{ }4,\mbox{ }5].$

The two common parametric solutions are \textit{(\ref{eq9})} and \textit{(\ref{eq10})}. The first is correct if $\frac{q}{p}=\frac{q_1 }{p_1 }$ is primary, 
and the second is correct if $\frac{q}{p}=\frac{q_2 }{p_2 }$ is secondary.
\begin{equation}
\label{eq9}
a=p^2 - q^2, \quad b=2pq,\quad c=p^2+q^2
\end{equation}
\begin{equation}
\label{eq10}
b=\frac{p^2 - q^2}{2},\quad a= \ pq,\quad c=\frac{p^2+q^2}{2}
\end{equation}
In our work we found the following mixed solution(s) quite convenient.
\begin{equation}
\label{eq11}
a=p_2 q_2, \quad b = 2p_1 q_1 ,\quad c=p_1 p_2 - q_1 q_2 = p_1 q_2 + p_2 q_1 
\end{equation}
For example, using key squence $[1,1,2,3]$: $$a = 3 =(1 \cdot 3),\quad b = 4 = 2(1\cdot 2), \quad c = 5 = (2 \cdot 3) - (1 \cdot 1) = (2 \cdot 1) + (3\cdot 1).$$ Write $G: \frac{q}{p} \to [a,b,c]$ to 
indicate that $\frac{q}{p}$ is a generator of $[a,b,c].$ Whether $\frac{q}{p}$ is primary or secondary will be stated, 
if not evident from the context.

The mixed form solution \textit{(\ref{eq11})} can be made more symmetric with 
the aid of the following definitions.
\[
r_1 = q_1 q_2 ,\quad r_2 = q_1 p_2,\quad r_3 = q_2 p_1,  \quad r_4 = p_1 p_2.
\]

\begin{flushleft}
\textbf{Circle Theorem:}
\end{flushleft}\emph{Triangle} $[a,b,c]$ \emph{has an in-circle with radius} $r_1$ \emph{and three ex-circles 
with radii} $r_2,r_3,r_4. \ $   \emph{Moreover,} (i) $r_1 + r_2 + r_3 = r_4,\quad r_1 \cdot r_4 = r_2 \cdot r_3.$ (ii) $a = r_1 + r_2 = r_4 -r_3,  \quad  b = r_1 + r_3 = r_4 - r_2,  \quad c = r_2 + r_3 = r_4 - r_1$
\begin{flushleft}
\mbox{ }
\end{flushleft}

It follows from \textit{(ii)} that four circles with radii $r_i$ form a 
tangent system with their centers at the corners of a rectangle of 
dimensions $a \times b.$ Proof of these claims is elementary. 
A more detailed correlation between tangent circles and ex/in-circles and 
full proofs can be found in [4]. See also Akhtar [1].

The key sequence $[3,1,4,5]$ was an example 
above. From it we can find the radii by multiplying out the product 
$(3+5)(1+4) \to  3,12,5,20.$ Then the radii give the sides of the 
triangle by  \textit {(ii)}: 

$a=3+12=20-5,$ \quad $b=3+5=20-12,$ \quad $c=12+5=20-3.$

The radii are useful here, primarily by means of the following. We assume 
$a<b.$

\begin{flushleft}
\textbf{Symphonic Corollary:}
\end{flushleft}
\begin{flushleft}
\textit{(i) The primary and secondary generators} $\rm T, T'$ \textit{for the symphonic major 
derivative are}
\[
{\rm T} = \tfrac{q_1 q_2 }{p_1 p_2 }= \tfrac{ab}{(c+a)(c+b)} = \tfrac{r_1 }{r_4 }
\]
\[
{\rm T'}= \tfrac{(r_2 + r_3 )}{(r_4 + r_1 )}= \tfrac{c}{(a+b)}.
\]
(\textit{ii) In the case of the minor derivative we have generators} $\rm T, T'$ \textit{given by}
\end{flushleft}
\[
{\rm T} = \tfrac{q_2 p_1 }{p_2 q_1 }=\tfrac{ab}{(c-a)(c+b)}= \tfrac{r_2 }{r_3},
\]
\[
{\rm T'}= \tfrac{(r_3 -r_2 )}{(r_3 + r_2 )} = \tfrac{(b-a)}{c}.
\]
We have put the most convenient expressions last on each line.

Remark: If $a>b,$ then in part\textit{(ii)} we must 
exchange $a, b$ and exchange $r_2 , r_3.$  Fairly routine computation 
using equations already given is all that is needed for proof.

Take as an example
\[
\rm S: \ [3,4,5] \to [35,12,37],	\quad S': \ [3,4,5] \to [5,12,13].
\]
Secondary and primary generators for $[3,4,5]$ are 
$\frac{1}{3},\mbox{ }\frac{1}{2}.$ Their product, $\frac{1}{6},$ and quotient,  
$\frac{2}{3},$ are the primary generators for the major and minor derivatives! 
Alternately, the radial quotients $\frac{r_1}{r_4},\mbox{ }\frac{r_2 }{r_3}$ give the same results.

We also check that $\frac{c}{a+b}=\frac{5}{7}$ and 
$\frac{b-a}{c}= \frac{1}{5}$ are the secondary generators. 
Note that $\frac{1}{6},\mbox{ }\frac{5}{7}$ combine in the 
key sequence $[5,\mbox{ }1,\mbox{ }6,\mbox{ }7],$ and similarly 
$\frac{2}{3}\mbox{ },\tfrac{1}{5}$ combine in $[1,\mbox{ 
}2,\mbox{ }3,\mbox{ }5].$

{\bf 3. A Symphonic Family History}

The set of PPT's has a family structure. It was independently 
discovered, seven years apart, by Barning [3] and Hall [7], so we refer to it as the {\bf {\it Barning-Hall 
Tree}}. It has often been rediscovered. The smallest triple $[3,4,5]$ is the \emph{root} and only occupant of level zero.  Level $n$ is converted to level $n+1$ by replacing each triple 
$[a,b,c]$ by its three immediate successors (its direct descendants, or ``children''). Hence level $n$ has $3^n$ members. Linking each triple with its three immediate successors creates a ternary, 
ordered, plane tree (compare \emph{Figure 4}).

Given the key sequence $[q_2,q_1,p_1,p_2]$ the three successors, called {\bf {\it left, middle, right}} are 
obtained by completing $[p_2,q_1,\cdots ], \ [p_2,p_1,\cdots ],$ $[q_2,p_1,\cdots].$ Thus our diagram agrees with Hall, and with Eckert [6] after it is rotated 90 degrees. 

Level {\bf {\it one}} is $[15,8,17], \quad [21,20,29], \quad [5,12,13].$ Level {\bf {\it 
two}} is 
\[
\begin{array}{*{20}r}
 {a:}   \hfill & {35} \hfill & {65} \hfill & {33} \hfill \\
 {b:}   \hfill & {12} \hfill & {72} \hfill & {56} \hfill \\
 {c:}   \hfill & {37} \hfill & {97} \hfill & {65} \hfill \\
\end{array} \begin{array}{*{20}r}
 \vert \hfill & {77} \hfill & {119} \hfill & {39} \hfill \\
 \vert \hfill & {36} \hfill & {120} \hfill & {80} \hfill \\
 \vert \hfill & {85} \hfill & {169} \hfill & {89} \hfill \\
\end{array} \begin{array}{*{20}r}
 \vert \hfill & {45} \hfill & {55} \hfill & { \ 7} \hfill \\
 \vert \hfill & {28} \hfill & {48} \hfill & {24} \hfill \\
 \vert \hfill & {53} \hfill & {73} \hfill & {25} \hfill \\
\end{array} 
\]
The (ternary) tree is simpler to display, and to analyze, when each 
triple is ``abbreviated'' to its primary generator: Let $\frac{q}{p}$ be any 
proper fraction. Consider the quotient $\frac{q}{p-2q}.$ This new quotient 
can be turned into a simpler positive proper fraction by (if necessary) 
changing sign, and taking the reciprocal, unless 
$\frac{q}{p}=\frac{1}{2}^{ }$or $\frac{q}{p}=\frac{1}{3}.$ 

Using this operation to define ``parent'' the proper fractions are formed 
into two ternary trees, one for primary generators and one for secondary 
generators. Only the first is needed or displayed here.

\begin{figure}[h]
\centerline{\includegraphics[width=4.70in,height=1.63in]{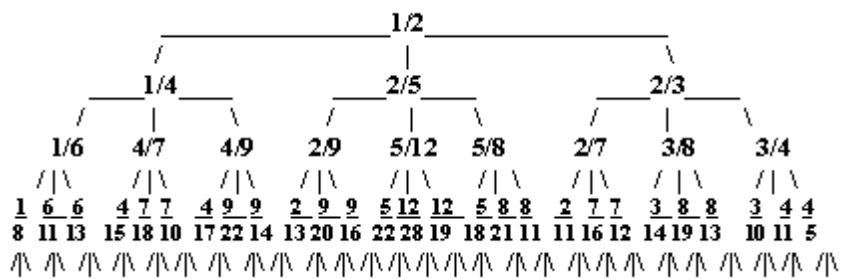}}
\label{fig4}
\caption{The Barning-Hall Tree (primary generators)}
\end{figure}

Any location on the tree can be reached by a path from the top, consisting 
of $n$ steps to reach the $n^{th}$ level. One step down-left (A), 
straight-down (B), or down-right (C), is given by 
\begin{equation}
\label{eq12}
\mbox{A: }\frac{q}{p}\mbox{ }\to \mbox{ }\frac{q}{p+2q},\mbox{ B: 
}\frac{q}{p}\mbox{ }\to \mbox{ }\frac{p}{2p+q},\mbox{ C: }\frac{q}{p}\to 
\mbox{ }\frac{p}{2p-q}.\mbox{ }
\end{equation}
Let us take the four triples from levels zero and one, and extract the major 
derivative. The primary generator is noted as a fraction.
\[
\mbox{ }\frac{1}{6}\mbox{ }[35,12,37];\mbox{ 
}\frac{3}{20}\mbox{ }[391,120,409];\mbox{ }\frac{6}{35}\mbox{ 
}[1189,420,1261];\mbox{ }\frac{2}{15}\mbox{ 
}[221,60,229]
\]
Likewise, the minor derivatives are as follows.
\[
\mbox{ }\frac{2}{3}\mbox{ }[5,12,13];\mbox{ 
}\frac{5}{12}\mbox{ }[119,120,169];\mbox{ }\frac{14}{15}\mbox{ 
}[29,420,421];\mbox{ }\frac{3}{10}\mbox{ 
}[91,60,109]
\]
The figure below shows the location of all eight of these on the tree. 
``Rho''($\rho$) marks the root. The major/minor derivatives are shown as small dark 
circles/rectangles (each circle is labeled with a generator).

\begin{figure}[h]
\centerline{\includegraphics[width=2.592in]{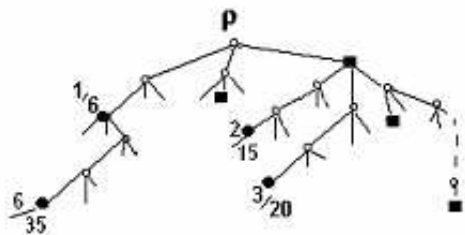}}
\label{fig5}
\caption{Locations of some major and minor derivatives.} 
\end{figure}


We can locate a fraction ($\frac{6}{35}$ say) on the tree by regression.
\[
\frac{6}{35}\mbox{ }\to \mbox{ }\frac{6}{23}\mbox{ }\to \mbox{ 
}\frac{6}{11}\mbox{ }\to \mbox{ }\frac{1}{6}\mbox{ }\to \mbox{ 
}\frac{1}{4}\mbox{ }\to \mbox{ }\frac{1}{2}\mbox{ }.
\]
Analyzing and reversing the steps gives us the code $\mbox{AACAA}$ for 
traveling ``down'' from the root, to arrive at the triple. The path codes 
for the major derivatives are $\mbox{AA, } \quad \mbox{CBAA, 
}\mbox{AACAA, } \quad \mbox{CAAA.}$ Those for the minor derivatives are 
$\mbox{C, BB, C}^{\mbox{13}}\mbox{, CCA.}$

By employing \textit{(\ref{eq12})}, these codes can be used as operators from the 
right side, e.g.
\[
\left( {\frac{\mbox{1}}{\mbox{4}}} \right)\mbox{CAA }= \mbox{ }\left( 
{\frac{\mbox{4}}{\mbox{7}}} \right)\mbox{AA }= \left( 
{\frac{\mbox{4}}{\mbox{15}}} \right)\mbox{A }= \left( 
{\frac{\mbox{4}}{\mbox{23}}} \right)\mbox{: [513, 184, 545].}
\]
Evidently most of the locations on the tree do not belong to major or minor 
derivatives.

The ternary tree is beautiful because it contains every primitive triple, 
each in its own unique position. It is possible to obtain trees using 
derivatives as well. On can define a binary tree in which the two immediate 
successors of a triple are the major derivative and the minor derivative. 
The chief reason we will not pursue this course is that there would be not 
one, but many such trees. Each triple that is not itself a derivative (major 
or minor) is the root of a different tree. 

It seems natural to pose certain questions:

\begin{enumerate}[(a)]
\item Can one characterize directly which triples are major or minor derivatives?
\item Can a triple be a derivative in two different ways?
\item For certain infinite sequences of triples, can we give the tree location of the major and minor derivatives?
\end{enumerate}
All in all, it seems best to start with the third question. We now add the 
symbols $\text S, \ \text{S}'$ to the one-step symbols $\mbox{A, B, C, }$ 
supplementing \textit{(\ref{eq12})} with \textit{(\ref{eq13})}.
\begin{equation}
\label{eq13}
{\rm S}: \tfrac{q}{p} \to \tfrac{q(p-q)}{p(p+q)}, \quad	{\rm S'}: \tfrac{q}{p}\to \tfrac{q(p+q)}{p(p-q)}
\end{equation}
In case of the latter expression, it is necessary to regard an improper fraction, say 
($\tfrac{3}{2}$), as an alternate generator $\frac{p}{q}$ 
equivalent to $\frac{q}{p}.$ (Thus $p = 3, \ q = 2$ either way).

This notation makes it convenient to investigate special cases,\ {\it viz.}
\[
\left( \tfrac{1}{2} \right){\rm A}^{n-1}{\rm S} =\left( \tfrac{1}{2n} 
\right){\rm S} = \tfrac{2n-1}{2n(2n+1)}.
\]
We have here a simple way of describing the symphonic major derivative of an 
arbitrary member of the family $[4n^2 - 1,\ 4n, \ 4n^2 + 1]$ attributed to Plato, and occupying the extreme left 
branch of the ternary tree. Likewise
\[
\left( \tfrac{1}{2} \right)\mbox{\rm C}^{n-1}\mbox{\rm S }=\left( \tfrac{n}{n+1} 
\right)\mbox{\rm S}=\tfrac{n}{(n+1)(2n+1)}
\]
describes the symphonic major derivative of the family on the extreme right 
branch of the ternary tree, attributed to Pythagoras. The Platonic family is 
characterized by $\mbox{A}^{n-1}\mbox{ }$or $\frac{1}{2n},$ and the 
Pythagorean family is characterized by ${\rm C}^{n-1},$ or 
$\frac{n}{n+1}.$ It is inevitable that we include also the ``Fermat'' 
family, characterized by ${\rm B}^{n-1}$ or $\frac{{\rm P}_n 
}{{\rm P}_{n+1} }$ (primary generator) or $\frac{{\rm Q}_n 
}{{\rm Q}_{n+1} }  $ (secondary generator) Here the Pell sequences 
\[
\begin{array}{l}
 (\mbox{P}_n)n=1,2,\cdots =\mbox{ }(1,\mbox{ }2,\mbox{ }5,\mbox{ 
}12,\mbox{ }29,\mbox{ }\cdots ,\mbox{ }a,\mbox{ }b,\mbox{ }2b+a,\mbox{ 
}\cdots ) \\ 
 (\mbox{Q}_n)n=1,2,\cdots =\mbox{ }(1,\mbox{ }3,\mbox{ }7,\mbox{ }17,\mbox{ 
}41,\mbox{ }\cdots ,\mbox{ }a,\mbox{ }b,\mbox{ }2b+a,\mbox{ }\cdots ) \\ 
 \end{array}
\]
are familiar recursive sequences, satisfying the indicated ``embedded'' 
recursion, which is not dissimilar to the Fibonacci recursion. We call 
attention to the key sequence formed from the generators: $$[\mbox{Q}_n 
,\mbox{ P}_n ,\mbox{ P}_{n+1} ,\mbox{ Q}_{n+1} ]\mbox{. }$$ Here the Fibonacci 
rule applies, and may be used to define both sequences. The ratios 
$\frac{\mbox{Q}_n }{\mbox{P}_n }$converge to $\sqrt 2 $; this helps to 
explain the ubiquity and popularity of these numbers. Because the Pell 
sequences are the subject of a large literature, we do not need to regret 
giving them such a brief glance.

The derived sequence$\mbox{ }\frac{1}{2}\mbox{P}_{2n} =\mbox{ Q}_n 
\mbox{P}_n $ was introduced above, just after the minor derivative, to 
explain the Fermat family in connection with Sastry's challenge problem! Working out the relationships (for instance Hatch [8]) we leave to the reader. The book by Barbeau [2] is 
very helpful, and thorough. 

We are now ready to continue, starting with the following.
\[
\begin{array}{l}
 \left( {\tfrac{1}{2}}\right){\text B}^{n-1}\ {\text S }\ = {\frac{\mbox{P}_n 
}{\mbox{P}_{n+1} }} \ {\text S }=\frac{\mbox{P}_{2n}     }{\mbox{P}_{2n+2} } 
\\ 
 \left({\tfrac{1}{2}}\right){\text B}^{n-1}\ {\text S'} = {\frac{\mbox{P}_n
}{\mbox{P}_{n+1} }} \ {\text S'}= \frac{\mbox{P}_{2n+1} \ - 1}
{\mbox{P}_{2n+1} \ + 1}\\

 \end{array}
\]
The final task is to specify the $\mbox{ABC}$ path (eliminate ${\text S, \ } {\text S}'$ from the path code). It is helpful to start by calculating a few cases 
and generalizing. The results are presented in the six tables of \emph{Figures 6} and \emph{7} below. Symbol $t$ 
is the primary generator, and $(t){\rm S}, \ (t){\rm S'}$ are the primary generators for the major/minor derivatives. 
Each has an associated path code.

One can verify from these tables that no duplication other than 
${\rm AA}$ occurs, forecasting the answer to 
question (b). It is also curious to observe that the Platonic family 
${\rm AA}\ldots $ is the trickiest, requiring a separation into even and 
odd cases.

\begin{figure}[h]
\centerline{\includegraphics[width=4.32in,height=6.0in]{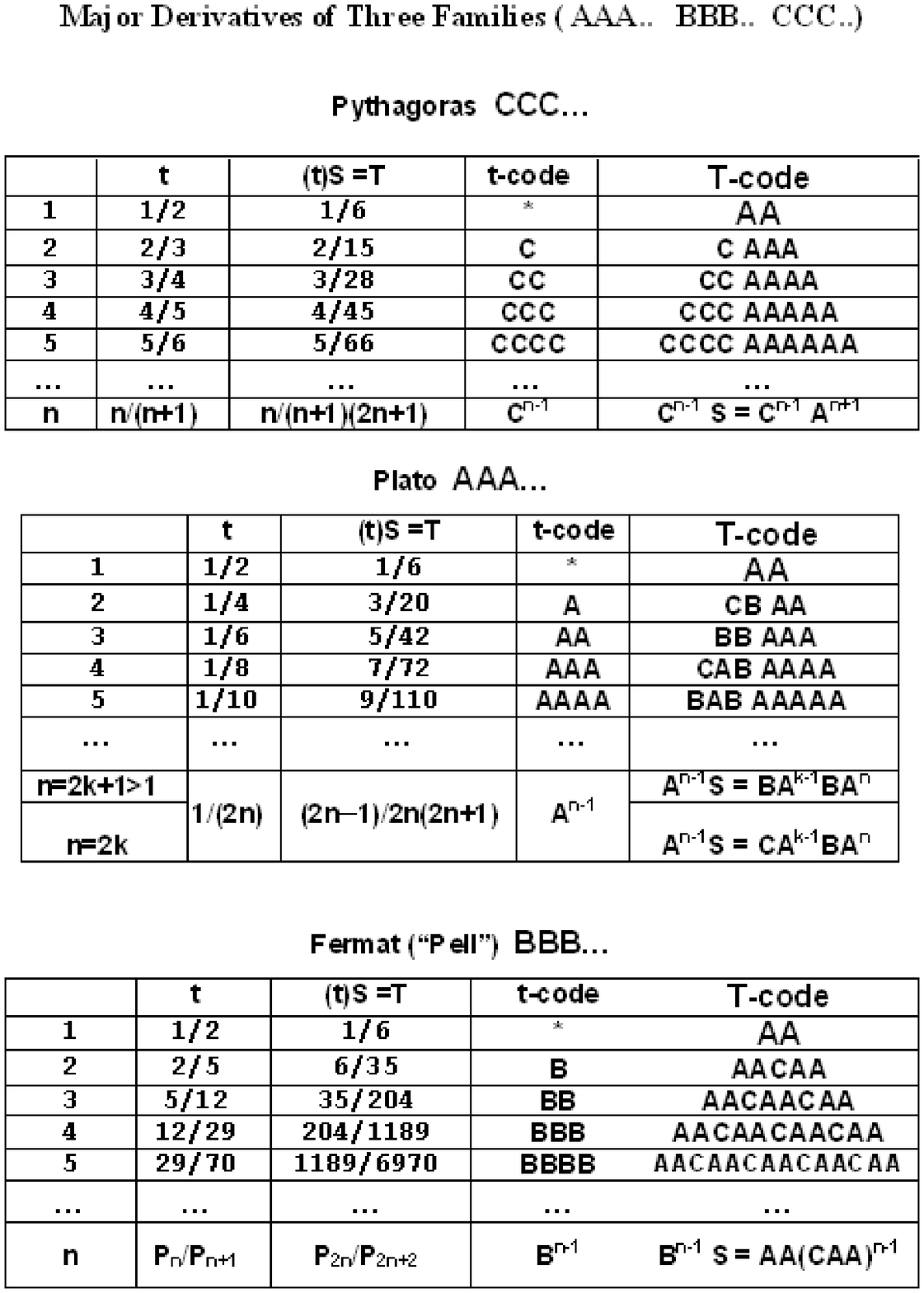}}
\label{fig6}
\caption{Three Tables: Major Derivative}
\end{figure}

\begin{figure}[h]
\centerline{\includegraphics[width=4.32in,height=6.42in]{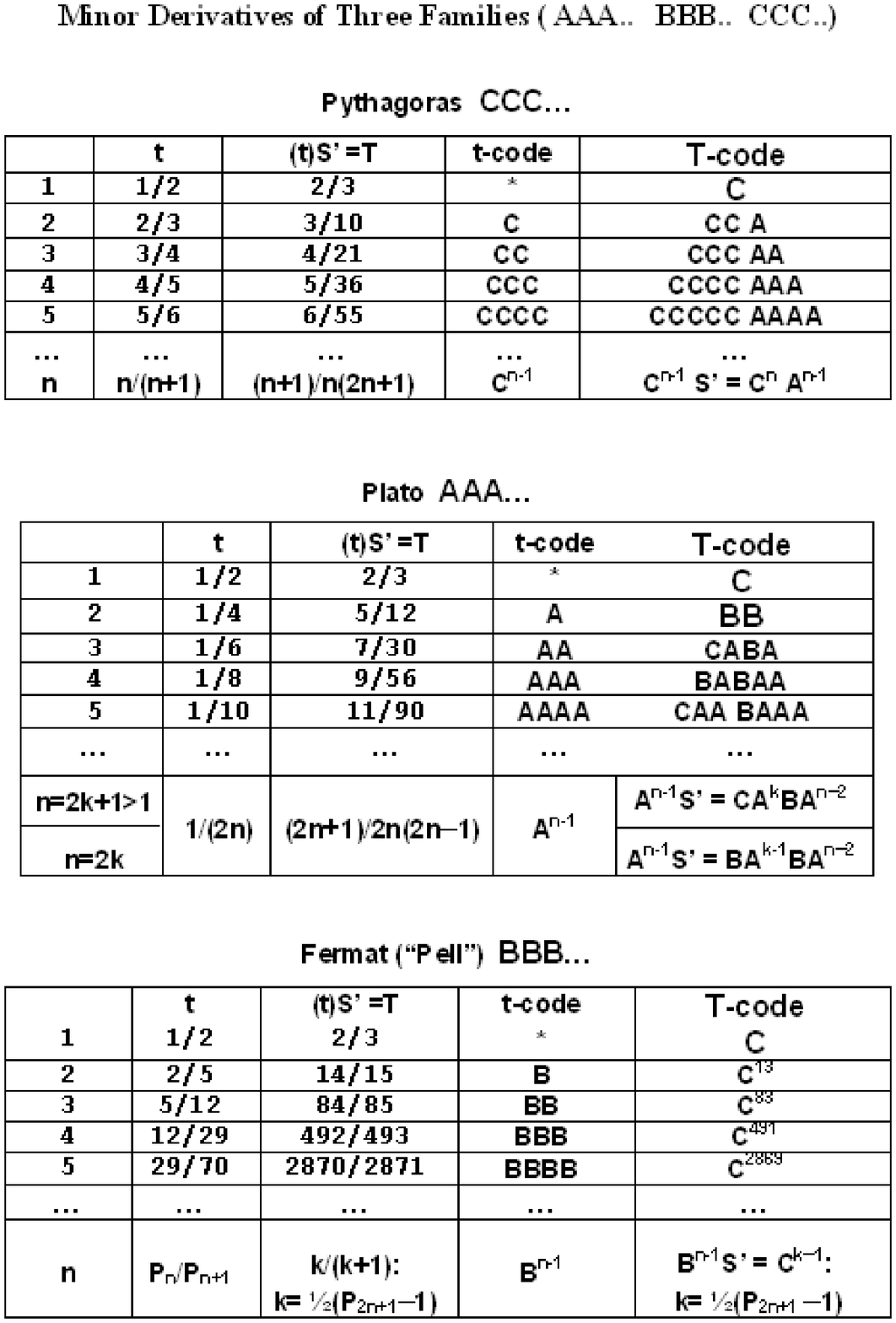}}
\label{fig7}
\caption{Three Tables: Minor Derivative}
\end{figure}

 A look at the final column of these tables shows that $\mbox{B}$ occurs only 
for the Platonic family and at most twice there!

When $\rm X,\ Y$ are path codes, write $\mbox{X}\cong \mbox{Y }$ to mean 
$( \tfrac{1}{2} )\rm X = (\tfrac{1}{2})\rm Y.$  If $\mbox{X, Y}$ are \emph{pure} codes, made up of $\mbox{A, B, C }$ only, then they are identical, letter for letter. But cases like $\rm AS'\cong \mbox{BB 
}$ are of interest here.

\begin{flushleft}
\textbf{Derivative Location Theorem:}
\end{flushleft}
 
\emph{The locations of the major and minor derivatives are given by the formula in the last row of the tables in Figures 6 and 7.}\newline \mbox{} \newline
We begin with the Pythagoras major derivative, claiming that 
$${\rm C}^{n-1} S \cong {\rm C}^{n-1} {\rm A}^{n+1}.$$ 
This is simple, for $\left( \frac{1}{2} \right){\rm C}^{n-1}= \frac{n}{n+1},$ 
and then $\left( {\tfrac{n}{n+1}} \right){\rm S} =\left( 
\tfrac{n}{n+1} \right)\left( \tfrac{1}{2n+1} \right), $ or 
alternately 
$$
\left( \tfrac{n}{n+1} \right)\mbox{A}^{n+1}= 
\tfrac{n}{\left({n+1} \right)+ 2n(n+1)}=\tfrac{n}{(2n+1)(n+1)}.
$$
The minor derivative ${\rm C}^{n-1} {\rm S'} \cong {\rm C}^n {\rm A}^{n-1}$ is similar. Note: direct 
descendents!

For the Fermat major derivative, recall Pell sequences $(\mbox{P}_n ),\mbox{ 
}(\mbox{Q}_n )$ and related sequence
\[
(\tfrac{1}{2})\mbox{P}_{2n} =(\mbox{P}_n \mbox{Q}_n )=(1,6,35,204,\ldots,a,b,c,\ldots),
\]
where $c = 6b-a$ is the recursion rule. First we have 
$\left( \tfrac {1} {2} \right){\rm B}^{n-1}=\tfrac{{\rm P}_n 
}{{\rm P}_{n+1} }.$ Here ${\rm T=}\tfrac{{\rm P}_n 
}{{\rm P}_{n+1} }$  is the primary generator, and $\rm T'=\tfrac{{\rm P}_{n+1} - {\rm P}_n }{{\rm P}_{n+1} + {\rm  P}_n }=\tfrac{ {{\rm Q}_n } }{ {{\rm Q}_{n+1} }}$ is the secondary 
generator. Then the major derivative is 
\[
\left( {\tfrac{{\rm P}_n }{{\rm P}_{n+1} }} \right){\rm S = TT' = }\left( \tfrac{{\rm P}_n }{{\rm P}_{n+1} } \right) \left( {\tfrac{{\rm Q}_n }{{\rm Q}_{n+1} }} \right) = \left( \tfrac { {\rm P_n}{\rm Q_n }} { {\rm P}_{n+1} {\rm Q}_{n+1} }\right) 
.
\]
Now $\left( \tfrac{1}{2} \right)\mbox{AA=}\left( \tfrac{1}{6} \right) = \tfrac{\rm P_1 Q_1} { \rm P_2 Q_2},$ and $\left( \tfrac{a}{b} 
\right)\mbox{CAA}=\tfrac{b}{6b-a}=\tfrac{b}{c},$ and so we get \newline
$\mbox{B}^{n-1} \mbox{S}\cong \mbox{AA(CAA)}^{n-1},$ since $\left( {\tfrac{{\rm P}_n }{{\rm P}_{n+1} }} \right){\rm CAA = }\left( 
{\tfrac{{\rm P}_{n+1} }{{\rm P}_{n+2} }} \right).$

The Fermat minor derivative requires $\mbox{B}^{n-1} {\rm S'}\cong 
\mbox{C}^{k-1}\mbox{.}$  Firstly, 
\[
\begin{array}{l}
 \left( \tfrac{ {\rm P}_n}{ {\rm P}_{n+1}} \right) {\rm S'} = \tfrac{\rm T }{\rm T'} = 
 \tfrac{ {\rm P}_n  {\rm Q}_{n+1}}   {{\rm P}_{n+1}{\rm Q}_n    } \quad n { \emph{  even}}        \\ 
 
 \left( \tfrac{ {\rm P}_n}{ {\rm P}_{n+1}} \right) {\rm S'} = \tfrac{\rm T'}{\rm T } = 
 \tfrac{{\rm P}_{n+1}{\rm Q}_n   }   { {\rm P}_n  {\rm Q}_{n+1} } \quad n  {\emph{ \  odd.}}       \\ 
 \end{array}
\] 
In both cases the final fraction has the form 
$\frac{k}{k+1}= \left( \frac{1}{2} \right)
\mbox{C}^{k-1}$ where 

$$
\begin{array}{rl}
k   = & \tfrac{1}{2} ({\rm P}_{2n+1} \ - 1 ) = min \{ {\rm P}_n \ {\rm Q}_{n+1},\quad {\rm P}_{n+1} \ \mbox{Q}_n \} \\
k+1 = & \tfrac{1}{2} ({\rm P}_{2n+1} \ + 1 ) = max \{ {\rm P}_n \ {\rm Q}_{n+1},\quad {\rm P}_{n+1} \ \mbox{Q}_n \} \\
\end{array}
$$

This is the special case of Sastry discussed earlier where triple $[a,b,c]$ with $b-a=1$ has minor 
derivative $[c,ab,ab+1].$ We also have $k=x+y$ where 
\[
(x,y)=(1,1),(6,8),(35,49),(204,288),\ldots 
\]
is any solution of the square-triangle problem: $x^2=y(y+1).$ The Platonic derivatives remain, and they require a distinction between odd and 
even. But the proof is relatively easy. To show $\mbox{A}_{n-1} 
\mbox{S}\cong \mbox{C A}^{k-1}\mbox{B A}^n\mbox{ }$ when $n=2k\mbox{ }$ is 
even, and $\mbox{A}^{n-1}\mbox{S}\cong \mbox{BA}^{k-1}\mbox{ BA}^n\mbox{ 
}$ when $n = 2k + 1 $ is odd, calculate
\[
\left( {\tfrac{1}{2}} \right)\mbox{A}^{n-1}\mbox{S }=\left( {\tfrac{1}{2n}} 
\right)\mbox{S}=\left( {\tfrac{1}{2n}} \right)\left( {\tfrac{2n-1}{2n+1}} 
\right)=\left( {\tfrac{2n-1}{2n(2n+1)}} \right),
\]
\[
\begin{array}{ll}
\left( {\tfrac{1}{2}} \right)\mbox{CA}^{k-1}= & \left( {\tfrac{2}{3}} 
\right)\mbox{A}^{k-1}=\left( {\tfrac{2}{3+4(k-1)}} \right)=\left( 
{\tfrac{2}{4k-1}} \right)=\left( {\tfrac{2}{2n-1}} \right)\mbox{ };\mbox{ 
}n\mbox{ }=\mbox{ }2k, \\

\left( {\tfrac{1}{2}} \right)\mbox{BA}^{k-1}= & \left( 
{\tfrac{2}{5}} \right)\mbox{A}^{k-1}=\left( {\tfrac{2}{5+4(k-1)}} 
\right)=\left( {\tfrac{2}{4k+1}} \right)=\left( {\tfrac{2}{2n-1}} 
\right)\mbox{ };\mbox{ }n\mbox{ }=\mbox{ }2k+1\mbox{,} \\
\end{array}
\]
$$
 \left( {\tfrac{2}{2n-1}} \right)\mbox{BA}^n=\left( {\tfrac{2n-1}{4n-2+\mbox{ 
}2}} \right)\mbox{A}^{n\mbox{ }}=\left( {\tfrac{2n-1}{4n}} 
\right)\mbox{A}^n= 
$$
$$
\left( {\tfrac{2n-1}{n(4n-2)+4n}} \right) 
 =\left( {\tfrac{2n-1}{2n(2n-1)+\mbox{ }2\mbox{ }}} \right)=\left( 
{\tfrac{2n-1}{2n(2n+1)}} \right).\mbox{ } 
$$

The case of $\mbox{A}^{n-1}{\rm S'}\cong 
\mbox{BA}^{k-1}\mbox{BA}^{n-2}, \quad n=2k;\mbox{ }$ 
$\mbox{A}^{n-1}\mbox{{S}'}\cong \mbox{CA}^k\mbox{BA}^n,\mbox{ }$ 
$n=2k+1\mbox{ }$is quite similar, and also routine. Q.E.D.

Taking a symphonic derivative at least doubles the length of the path code, 
which supports the empirical conclusion that these derivatives are quite sparse. 
The Fermat minor derivatives are especially far down the tree!

\begin{flushleft}
{\bf 4. Musical Aptitude}
\end{flushleft}

Just how improbable is it for an arbitrary triple to be symphonic? This is 
another approach to question (a) above. The following development 
shows that often it is easy to conclude that a triple is not symphonic.

For any PPT $[a,b,c] $ the product $abc$ is 
divisible by $60.$ This fact is \emph{long} known (Sierpinski [11]) but apparently not always \emph{well} known (Monaghan [9]). 
Recall that we assume that $a$ is odd. Write $x\vert y$ to say that 
integer $x$ divides integer $y.$ Then
\begin{equation}
\label{eq14}
4\vert b, \quad	3\vert ab,\quad	5\vert abc.
\end{equation}
By the mixed solution \textit{(\ref{eq11})}, we can write
\[
b=2q_1 p_1, \quad ab=2q_1 p_1 q_2 p_2, \quad abc=2q_1 p_1 q_2 p_2 (q_1 p_2 + p_1 q_2 ).
\]

Let ${\rm Q}=(\ldots ,x,y,z,u,v,\ldots )$ be a Fibonacci rule 
sequence. Say it is $k\mbox{-}$primitive if not every term is divisible by 
$k.$ Then if $\rm Q$ is $\mbox{2-}$primitive, $2\vert xyz$ and every {\it third} item is even. If $\rm Q$ is $\mbox{3-}$primitive, $3\vert 
xyzu$ and every {\it fourth} item is a multiple of three, and if $\rm Q$ is 
$\mbox{5-}$primitive, either $5\vert xyzuv$ and every {\it fifth} item is a 
multiple of five, or else {\it no} item is a multiple of five. The last possibility 
must equal $(\ldots ,\ 1,-2,-1, \ 2,\ 1,\ldots )$ mod(5).

Since key sequence $(q_2,q_1,p_1,p_2 )$ has the Fibonacci property, and $q_2 $ is odd, it is easy to show that 
$b$ is a multiple of four. Also $ab$ is divisible by three (and $c$ therefore cannot be). Finally $5\vert abc,$ -- the trickiest case.

The key sequence may be extended indefinitely at either end, yielding a 
longer sequence with the Fibonacci property. If $5\vert ab,$ we are 
done. Else (i) the sequence is a part of a longer sequence $(5x,q_2 
,q_1,p_1,p_2,5y)$ and is equivalent $\bmod (5)$ to $(0,x,x,2x,3x,0),$ or (ii) extends 
indefinitely without including a multiple of five.

Suppose case (i). Given that $c=(q_1 p_2 + p_1 q_2 ),$ this is $\bmod (5)$ congruent 
to $(x)(2x) + (x)(3x)= x^2(2+3)= 0.$ Now suppose  
(ii). The key sequence is a four term subsequence of $(\ldots , \ 1,-2,-1, \ 2, \ 1,\ldots ).$ Wherever the starting point, in $\bmod (5)$ we get $c=(2)(-2)+(1)(-1)=0.$ Q.E.D.

\begin{table}[htbp]
\begin{center}
Based on the location of the factors, we can list six classes of PPT's: \newline
$\mbox{}$ \newline
\begin{tabular}{|l|l|l|}
\hline
 & 
{\bf 3$\vert $a}& 
{\bf 3$\vert $b} \\
\hline
{\bf 5$\vert $c }& 
{\bf T1: [ 3x, 4x, 5x] }& 
{\bf T2: [ x, 12y, 5z]} \\
\hline
{\bf 5$\vert $a}& 
{\bf T3: [15x, 4x, x] }& 
{\bf T4: [5x, 12y, z]} \\
\hline
{\bf 5$\vert $b}& 
{\bf T5: [ 3x, 20y, z] }& 
{\bf T6: [ x, 60y, z]} \\
\hline
\end{tabular}
\label{tab1}
\end{center}
\end{table}

\begin{table}[h]
\begin{center}
The simplest member of each (infinite!) class is:\newline
$\mbox{}$ \newline
\begin{tabular}{|l|l|l|l|l|l|}
\hline
$^{1}${\bf /}$_{2}$& 
$^{3}${\bf /}$_{4}$& 
$^{1}${\bf /}$_{4}$& 
$^{2}${\bf /}$_{3}$& 
$^{2}${\bf /}$_{5}$& 
$^{5}${\bf /}$_{6}$ \\
\hline
{\bf [3,4,5]}& 
{\bf [7,24,25]}& 
{\bf [15,8,17]}& 
{\bf [5,12,13]}& 
{\bf [21,20,29]}& 
{\bf [11,60,61]} \\
\hline
\end{tabular}
\label{tab2}
\end{center}
\end{table}

A symphonic derivative can only be in class $\mbox{T4 }$or $\mbox{T6.}$ More 
specifically, using 
$
{\rm S, S':} [a,b,c] \to [\rm A, B, C] = [\vert a\pm b\vert c, \ ab, \ c^2 \pm ab]
$
we see that factors of $c$ are transferred to $\rm A,$ and factors of 
$a,b$ are transferred to $\rm B.$ 

\begin{flushleft}
\textbf{Symphonic Factor Theorem:}
\end{flushleft}
\emph{The symphonic product} $[\rm A, B, C]$ \emph{is in 
class} $\rm T4$ \emph{if the original} $[a,b,c]$ \emph{is in} $\rm T1$ or $\rm T2,$ \emph{but otherwise 
is in class} $\rm T6.$ \emph{Hence after two steps, a symphonic product is in} $\rm T6.$ \newline \mbox{ }

Are all class $\rm T6$ PPT's symphonic derivatives?  We explore this question making use of a famous PPT studied by Fermat.  According to Sierpinski [11], Fermat wrote a letter to Mersenne in the 
year 1643, in which he asserted that the PPT
\[
[\mbox{ }456\mbox{ }54860\mbox{ }27761,\mbox{ }106\mbox{ }16522\mbox{ 
}93520,\mbox{ }468\mbox{ }72986\mbox{ }10289\mbox{ }]
\]
is the \emph{smallest} PPT in which the hypotenuse and the sum of the legs are both 
squares. Just for starters, we compute the primary generator and use it to 
locate the tree position of Fermat's triple.

The generator turns out to be $\raise0.7ex\hbox{${\mbox{246792}}$} 
\!\mathord{\left/ {\vphantom {{\mbox{246792}} 
{\mbox{2150905}}}}\right.\kern-\nulldelimiterspace}\!\lower0.7ex\hbox{${\mbox{2150905}}$}.$ 
Now we back up through the tree $41$ steps! (see Table 1)

\begin{table}[h]
\begin{center}
\textbf{The 41 Steps}

\begin{tabular}{|p{175pt}|l|}
\hline
{ A \quad 246792 / 1657321 \quad (\emph{Fermat})} &  
{ C \quad 3755 / 5778}        \quad \quad \quad \quad \quad \quad \mbox { }\\
\hline
{ A \quad 246792 / 1163737 }& 
{ C \quad 732 / 3755 } \\
\hline
{ A \quad 246792 / 670153 }& 
{ B \quad291 / 1732} \\
\hline
{ B \quad 176569 / 246792 }& 
{ A \quad 291 / 1150 } \\
\hline
{ C \quad 106346 / 176569 }& 
{ A \quad 291 / 568}  \\
\hline
{ C \quad 36123 / 106346} & 
{ C \quad 14 / 291} \\
\hline
{ B \quad 34100 / 36123 }& 
{ A \quad 14 / 263 } \\
\hline
{ C \quad 32077 / 34100 }& 
{ A \quad 14 / 235 } \\
\hline
{ C \quad 30054 / 32077 }& 
{ A \quad 14 / 207} \\
\hline
{ C \quad 28031 / 30054 }& 
{ A \quad 14 / 179 } \\
\hline
{ C \quad 26008 / 28031 }& 
{ A \quad 14 / 151 } \\
\hline
{ C \quad 23985 / 26008 }& 
{ A \quad 14 / 123} \\
\hline
{ C \quad 21962 / 23985 }& 
{ A \quad 14 / 95 } \\
\hline
{ C \quad 19939 / 21962 }& 
{ A \quad 14 / 67 } \\
\hline
{ C \quad 17916 / 19939 }& 
{ A \quad 14 / 39} \\
\hline
{ C \quad 15893 / 17916 }& 
{ B \quad 11 / 14 } \\
\hline
{ C \quad 13870 / 15893 }& 
{ C \quad 8 / 11 } \\
\hline
{ C \quad 11847 / 13870 }& 
{ C \quad 5 / 8} \\
\hline
{ C \quad 9824 / 11847 } & 
{ C \quad 2 / 5} \\
\hline
{ C \quad 7801 / 9824 }  & 
{ B \quad 1 / 2 \quad (\emph{Root})} \\
\hline
{ C \quad 5778 / 7801 }& 
{ ....\quad q/p} \\
\hline 
\end{tabular}
\label{}
\caption{Path code to Fermat's triangle}
\end{center}
\end{table}


Reversing and collating the sequence of letters gives the path code locator 
for Fermat's enormous triple. It has an interesting structure that may 
perhaps be typical.
\[
{\underbrace{\rm BCCCB}_{5} \underbrace{\rm AAAAAAAAA}_{9} \underbrace{\rm CAAB}_{4} \underbrace{\rm CCCCCCCCCCCCCCCC}_{16} \underbrace{\rm BCCB}_{4} \underbrace{\rm AAA}_{3}}
\]
\[
\mbox{ 41= 5 + 9 + 4 + 16 + 4 + 3}
\]
Since Fermat's triple falls in class $\mbox{T6}$ we ask if it is a symphonic 
derivative. First we factor the numbers (from here on, A,B,C represent numbers):
\[
\begin{array}{rl}
{\rm A   =}& 456 54860 27761 = 17 \cdot 31 \cdot 239 \cdot 257 \cdot 141041 \\ 
{\rm B   =}& 106 16522 93520 = 16 \cdot 3 \cdot 5 \cdot  7 \cdot 13 \cdot 113 \cdot 277 \cdot 1553 \\ 
{\rm C   =} & 468 72986 10289 = (2165017)^{2} \\ 
{\rm A+B =} & 562 71383 21281 = (1009 \cdot 2351)^{2} \\ 
{\rm C B =} & 362 56463 16769 = (31 \cdot 239 \cdot 257)^{2} \\ 
 \end{array}
\]
Assume, for arguments sake, that $[\rm A, B, C] =[(a+b)c, \ ab, \ ab + c^2].$ This then 
implies that $\rm CB =c^2$ and $a+b = 17 \cdot 141041.$ Thus 
\[
(a+b)^2-4ab = (17 \cdot 141041)^2 -4{\rm B }= 313 \cdot 383 \cdot 12532151.
\]
This can not be the square $(a-b)^2$ thus Fermat's triple is 
not ``major''. Alternately we find that $c=31 \cdot 239 \cdot 257,$ which contains isolated 
primes of form $4k+3,$ contrary to known results. With the necessary 
adjustments, a similar argument establishes that the Fermat triple is not a 
minor derivative as well.

Finally, we conclude with some results on query (b).

\begin{flushleft}
\textbf{Anti-derivative Theorem}
\end{flushleft}
\emph{If }$[a,b,c] \to [\rm A, B, C]$ \emph{by either} 
$\rm S$ or $\rm S'$ then $[\rm A, B, C]$ \emph{determines} $[a,b,c].$ \newline \mbox{  }

Let $\frac{\rm Q}{\rm P}$ be the primary generator of the 
derivative $[\rm A, B, C].$ Suppose that
\[
[\rm A, B, C] = [{\rm P}^2 - {\rm Q}^2, \ 2PQ, \ {\rm P}^2 + {\rm Q}^2] = [c(a+b),\ ab, \ c^2 + ab].
\]
 $${\rm Then:}\quad {\rm P-Q} = c, \quad  {\rm P+Q} = a+b, \quad  {\rm 2PQ} = ab.$$ The 
quadratic $x^2  - ({\rm P} + {\rm Q})x + {\rm 2PQ} $ has roots $a,b.$ Now suppose
\[
[\rm A, B, C] = [{\rm P}^2 - {\rm Q}^2,\ 2PQ, \ {P}^2 + {Q}^2] = [c(b-a),\ ab,\ c^2 - ab].
\]
 $$ {\rm Then:} \quad {\rm P+Q} = c,\quad  {\rm P-Q} = b-a,\quad {\rm 2PQ} = ab,
$$
(We have assumed $a<b$ which could have a side affect -- the exchange 
of radii $r_2, r_3.$)  From quadratic $x^2 + ({\rm P} - {\rm Q}) x - {\rm 2PQ }$ we find $(b,-a).$

This result allows us to compute ``anti-derivatives'' for \emph{any} primitive triple! The triple ${\rm T} = [15, 8, 17]$ is not a major or a minor derivative. Undaunted, we produce real and complex surds as 
the anti-derivatives.
\[
\begin{array}{lrl}
 {\rm S }^{-1}(\rm T) = [a,b,3], \quad &  a,b = & \tfrac{1}{2}  \left(5 \pm \sqrt {-7} \right), \\ 
 {\rm S'}^{-1}(\rm T) = [a,b,5], \quad & -a,b = & \tfrac{1}{2}  \left(3 \pm \sqrt {\ 41} \right). \\ 
 \end{array}
\]
These are ``right triangles'' that transform by ${\rm S, S'}$ to 
$[15, 8, 17].$

\begin{flushleft}
\textbf{Major/Minor Theorem:}
\end{flushleft}
 \emph{There is no triple} 
$[\rm A, B, C]$ \emph{which is both major and minor.} \newline \mbox{ }

Proceeding in a similar manner, suppose $[a,b,c]$ is a major 
anti-derivative, and suppose $[a',b',c']$ is a minor 
anti-derivative for $[\rm A, B, C].$ Then
\[
 {\rm P + Q} = c' = a+b, \quad {\rm P - Q} = c = b'-a', \quad 2{\rm PQ }= ab = a'b'.
\]
The last pair of equations imply that the two anti-derivatives are triangles 
with the same area. The first two equations can be added and subtracted to 
show that $\{r_1,r_4 \} = \{r'_2, r'_3\}$ which with the equations 
$r_1 r_4 = r_2 r_3, \quad r'_1 r'_4 = r'_2 r'_3$ implies $\{r'_1, r'_4\} = \{r_2, r_3\}.$ But this is impossible, since $r_2,r_3 < r_4$ and $r'_2, r'_3 < r'_4.$

\newpage

\newpage
\begin{center}
{\bf Author Contact Information}
\end{center}

\small{H. Lee Price, 83 Wheatstone Circle, Fairport, NY 14450-1138.

email: \textit{tanutuva@rochester.rr.com} }\newline

\small{Frank R Bernhart, Math Dept Visitor, University of Illinois, Urbana, IL 61801.

email: \textit{bernhart@math.uiuc.edu}}

\end{document}